\documentclass[reqno,11pt]{amsart}

\usepackage{amsmath,amsfonts,amssymb,amsthm,epsfig,amscd,}
\usepackage[]{fontenc}
\usepackage[latin1]{inputenc}
\usepackage{enumerate}
\usepackage[cmtip,all]{xy}
\usepackage{tabularx,multirow}

 \allowdisplaybreaks \numberwithin{equation}{section}

\newtheorem{theorem}{Theorem}[section]

\newtheorem{lemma}[theorem]{Lemma}

\theoremstyle{definition}
\newtheorem{definition}[theorem]{Definition}

\theoremstyle{remark}
\newtheorem{remark}[theorem]{Remark}
\newtheorem{example}[theorem]{Example}

\DeclareMathOperator{\irr}{irr}
\DeclareMathOperator{\stabirr}{stab.irr}
\DeclareMathOperator{\unirr}{uni.irr}
\DeclareMathOperator{\conngon}{conn.gon}
\DeclareMathOperator{\covgon}{cov.gon}
\DeclareMathOperator{\gon}{gon}

\begin{document}

\title[On irrationality of surfaces in $\mathbb{P}^3$]{On irrationality of surfaces in $\mathbb{P}^3$}
\author{Francesco Bastianelli}
\address{Dipartimento di Matematica, Universit\`a degli Studi di Bari, Via Edoardo Orabona 4, 70125 Bari -- Italy}
\email{francesco.bastianelli@uniba.it}
\thanks{This work was partially supported by FIRB 2012 \emph{``Spazi di moduli e applicazioni''}; MIUR PRIN 2010--2011  \emph{``Geometria delle variet\`a algebriche''}; INdAM (GNSAGA)}

\begin{abstract}
The degree of irrationality $\irr(X)$ of a $n$-dimensional complex projective variety $X$ is the least degree of a dominant rational map $X\dashrightarrow \mathbb{P}^n$.
It is a well-known fact that given a product $X\times \mathbb{P}^m$ or a $n$-dimensional variety $Y$ dominating $X$, their degrees of irrationality may be smaller than the degree of irrationality of $X$. 
In this paper, we focus on smooth surfaces $S\subset\mathbb{P}^3$ of degree $d\geq 5$, and we prove that $\irr(S\times\mathbb{P}^{m})=\irr(S)$ for any integer $m\geq 0$, whereas $\irr(Y)<\irr(S)$ occurs for some $Y$ dominating $S$ if and only if $S$ contains a rational curve.
\end{abstract}

\maketitle

\section{Introduction}\label{section INTRODUCTION}

In the recent paper \cite{BDELU}, several perspectives for studying measures of irrationality for projective varieties have been proposed.
Along the same lines, we discuss various birational invariants, which extend the notion of gonality to higher dimensional varieties and, more importantly, they somehow measure the failure of a given variety to satisfy certain rationality properties.
In particular, we focus on smooth surfaces in $\mathbb{P}^3$ of degree $d\geq 5$, and we complete the characterization of the invariants, depending on the degree $d$ and on the existence of special subvarieties. 

Given a smooth complex projective curve $C$, the \emph{gonality} of $C$ is defined as the least integer $\delta$ such that there exists a non-constant morphism $C\longrightarrow \mathbb{P}^1$ of degree $\delta$, and it is denoted by $\gon(C)$. 
This is one of the most studied and important invariants in the theory of algebraic curves, and it can be thought as measuring how far the curve $C$ is from being rational.
The most natural extension of the notion of gonality to higher dimensional varieties is probably the \emph{degree of irrationality}.
For a smooth complex projective variety $X$ of dimension $n$, it is defined as
\begin{displaymath}
\irr(X):=\min\Big\{\delta\in \mathbb{N}\Big|
\exists  \text{ a dominant rational map }X\dashrightarrow \mathbb{P}^n \text{ of degree }\delta
\Big\}.
\end{displaymath}
Initially, the degree of irrationality was introduced in terms of field extensions by Heinzer and Moh, in order to discuss some generalizations of L\"uroth Theorem (cf. \cite{MH1,MH2}).
It is also worth mentioning a series of papers by Yoshihara, where this invariant was studied especially in the case of surfaces (see e.g. \cite{TY,Y1,Y2,Y3,Y4}).
We note further that the degree of irrationality equals 1 if and only if $X$ is a \emph{rational} variety.

In the case of smooth surfaces in $\mathbb{P}^3$, the degree of irrationality is governed by the following (see \cite[Theorem 1.3]{BCD}).
\begin{theorem}\label{theorem Cortini}
Let $S\subset \mathbb{P}^3$ be a smooth surface of degree $d\geq 5$.
Then $\irr(S)=d-2$ if and only if one of the following occurs
\begin{itemize}
  \item[(a)] $S$ contains a twisted cubic; 
  \item[(b)] $S$ contains a rational curve $R$ of degree $r$ and a line $\ell$ which is $(r-1)$-secant to $R$.
\end{itemize}
Otherwise, $\irr(S)=d-1$.
\end{theorem}
Moreover, when $S\subset \mathbb{P}^3$ is a very general surface of degree $d\geq 6$, then $\irr(S)=d-1$, and it is computed only by projections from points $p\in S$. This fact generalizes a famous result due to M. Noether on the gonality of plane curves (cf. \cite{C,H}), and for any $n\geq 3$, it has been recently extended to very general hypersurfaces $X\subset \mathbb{P}^{n+1}$ of degree $d\geq 2n+2$ (see \cite{BDELU}).

\smallskip
It is worth noting that the degree of irrationality may decrease when we consider a product $X\times \mathbb{P}^{m}$, or a $n$-dimensional variety $Y$ dominating $X$.
Namely, it can happen that $\irr(X\times\mathbb{P}^{m})<\irr(X)$ and $\irr(Y)<\irr(X)$.
It is indeed a remarkable established fact the existence of stably rational and unirational varieties which are not rational (see e.g. \cite{AM,BCSS,CG,IM}).   
In order to investigate when these phenomena occur, we consider the invariants
\begin{align*}
& \stabirr(X)  :=\min\Big\{\irr(X\times \mathbb{P}^m)\Big|
m\in \mathbb{N}\Big\}\\
& \unirr(X)  :=\min\left\{\irr(Y)\left|\begin{array}{l}
\exists  \text{ a dominant rational map }Y\dashrightarrow X \\
\text{with }\dim Y=\dim X\end{array}
\right.\right\},
\end{align*}
and we compare them with $\irr(X)$. 
We note that the conditions $\unirr(X)=1$ and  $\stabirr(X)=1$ recover the definitions of \emph{stably rational} varieties and \emph{unirational} varieties, respectively.
Furthermore, it is easy to check that $\irr(X)\geq \stabirr(X)\geq \unirr(X)$  (see Lemma \ref{lemma Inequalities}).

We characterize these invariants on smooth surface in $\mathbb{P}^3$, and we achieve the following.
\begin{theorem}\label{theorem Main}
Let $S\subset \mathbb{P}^3$ be a smooth surface of degree $d\geq 5$.
Then 
\begin{itemize}
  \item[(i)] $\displaystyle\stabirr(S)=\irr(S)$.
  
  \smallskip
  \item[(ii)] $\unirr(S)=d-2$ if and only if $S$ contains a rational curve. Otherwise, $\unirr(S)=d-1$.
\end{itemize}
\end{theorem}
\noindent In particular, it follows that $\irr(X\times \mathbb{P}^{m})$ cannot drop, and $\stabirr(S)$ is governed by Theorem \ref{theorem Cortini}.
Furthermore, although unirationality is equivalent to rationality in the case of surfaces, the equality $\unirr(S)=\irr(S)$ may fail when $\unirr(S)>1$.
In particular, Theorems \ref{theorem Cortini} and \ref{theorem Main} describe when this phenomenon occurs for smooth surfaces in $\mathbb{P}^3$ (cf. also \cite[Remark 4.10]{BCD}). 

\smallskip
To conclude our analysis, we would like to view our result in a more general setting.
So we consider two further invariants, which were introduced in \cite{BDELU}; the \emph{covering gonality}
\begin{displaymath}
\covgon(X):=\min\left\{c\in \mathbb{N}\left|
\begin{array}{l}
\text{Given a general point }x\in X,\,\exists\text{ an irreducible}\\ \text{curve } C\subset X  \text{ such that }x\in C \text{ and }\gon(C)=c
\end{array}\right.\right\},
\end{displaymath}
and the \emph{connecting gonality}
\begin{displaymath}
\conngon(X):=\min\left\{c\in \mathbb{N}\left|
\begin{array}{l}
\text{Given two general points }x,y\in X,\,\exists \text{ an irreducible}\\  \text{curve } C\subset X \text{ such that }x,y\in C \text{ and }\gon(C)=c
\end{array}\right.\right\}.
\end{displaymath}
Here the curves involved in the definitions are possibly singular, so $\gon(C)$ stands for the gonality of the normalization of $C$.
Moreover, the conditions $\covgon(X)=1$ and $\conngon(X)=1$ characterize \emph{uniruled} varieties and \emph{rationally connected} varieties, respectively.

We point out that all the five notions we introduced above are birational invariants, and when $X$ is a curve, each of them does coincide with gonality (cf. Remark \ref{remark Curves}).
Furthermore, any such an invariant may be thought as a measure of the failure of $X$ to satisfy the corresponding rationality property appearing in the classical chain of implications
\begin{equation*}
\text{rational} \implies \text{stably rational} \implies\text{unirational} \implies\text{rationally connected} \implies\text{uniruled}.
\end{equation*}
Accordingly, these `measures of irrationality' fit in the sequence of inequalities
\begin{equation*}
\irr(X)\geq \stabirr(X)\geq\unirr(X)\geq\conngon(X) \geq\covgon(X).
\end{equation*}
Of course, several classical results can be rephrased in terms of these invariants (see e.g. \cite[p. 351]{MM} and \cite[p. 833]{D} in the case of K3 surfaces), and other classes of varieties can be easily handled, as in Example \ref{example Ruled Surfaces}.
However, the problem of determining these invariants is in general widely open, even in the case of surfaces.
We suggest for instance \cite{B2} for a survey on the topic, and \cite[Section 4]{BDELU} for a series of very interesting problems on various measures of irrationality for projective varieties.

Turning to smooth surfaces $S\subset \mathbb{P}^3$ of degree $d$, we have that when $1\leq d\leq 3$, the surfaces are rational, so the problem is trivial as all the invariants equal 1.
If instead $d=4$, the invariants are completely characterized, except for $\stabirr(S)$, but this is actually the unique case we can not decide (see Remark \ref{remark K3}). 
On the other hand, when $S\subset \mathbb{P}^3$ is a smooth surface of degree $d\geq 5$, Lopez and Pirola proved that the covering gonality is $\covgon(S)=d-2$, and they classified all possible families of $(d-2)$-gonal curves covering $S$ (see \cite[Corollary 1.7]{LP}). 
Apart from those families depending on the existence of rational and elliptic curves on $S$, every such a surface is covered by the family of $(d-2)$-gonal curves obtained as tangent hyperplane sections, $C_p:= S\cap T_pS$ with $p\in S$.
In fact, it is easy to check that this family computes also the connecting gonality, that is $\conngon(S)= d-2$ (cf. \cite[Example 1.7]{BDELU}).
Therefore, taking into account Theorems \ref{theorem Cortini} and \ref{theorem Main}, the problem is completely understood, and the picture is summarized by the following.
\begin{theorem}
Let $S\subset \mathbb{P}^3$ be a smooth surface of degree $d\geq 5$.
Then 
\begin{itemize}
  
  \smallskip
  \item[(i)] $\covgon(S)=\conngon(S)=d-2$,
  
  \smallskip
  \item[(ii)] $\displaystyle\unirr(S)=\left\{\begin{array}{ll} d-2 & \text{if $S$ contains a rational curve}\\ d-1 & \text{otherwise,}\end{array}\right.$
  
  \smallskip
  \item[(iii)] $\displaystyle\stabirr(S)=\irr(S)=\left\{\begin{array}{ll} d-2 & \text{if either }\mathrm{(a)}\text{ or }\mathrm{(b)}\text{ in Theorem \ref{theorem Cortini} occurs}\\ d-1 & \text{otherwise.}\end{array}\right.$
\end{itemize}
\end{theorem}

In particular, when $S$ is assumed to be very general, then $\covgon(S)=\conngon(S)=d-2$ and $\unirr(S)=\conngon(S)=\irr(S)=d-1$, as $S$ does not contain rational curves (see e.g. \cite{X}).
Besides, in the light of \cite[Theorem A]{BDELU} and \cite[Theorem 3.3]{BCFS}, it would be interesting to understand the behavior of those invariants when $X\subset \mathbb{P}^{n+1}$ is a very general hypersurface of large degree and arbitrary dimension.

The proof of Theorem \ref{theorem Main} relies mainly on the classification of correspondences with null trace on smooth surfaces in $\mathbb{P}^3$, which is described by \cite[Theorem 1.3]{LP} and is based on Mumford's technique of induced differentials (see \cite[Section 2]{M}).  
In the next section, we thus discuss some properties of the birational invariants we introduced, and we follow \cite{LP} in order to relate them to correspondences with null trace. 
Then, Section 3 shall be entirely devoted to prove Theorem \ref{theorem Main}.

\section{Preliminaries}

\subsection*{Notation}
We shall work throughout over the field $\mathbb{C}$ of complex numbers.
By \emph{variety} we mean a complete reduced algebraic variety over $\mathbb{C}$, unless otherwise stated.
We say that a property holds for a \emph{general} point ${x\in X}$ if it holds on an open non-empty subset of $X$.
Analogously, we say that a property is satisfied by a \emph{very general} point ${x\in X}$, if the locus of points sharing the property is the complement of a countable collection of proper subvarieties of $X$.

\subsection{Measures of irrationality}\label{subsection Invariants}

In this subsection, we discuss some elementary properties of the invariants we defined in the Introduction.

Let $X$ be a smooth projective variety of dimension $n$.
By composing birational maps $X'\dashrightarrow X$ and dominant rational maps $X\dashrightarrow \mathbb{P}^n$, it is trivial to check that the degree of irrationality $\irr(X)$ is a birational invariant, and analogously, also $\stabirr(X)$ and $\unirr(X)$ are. 

On the other hand, the covering gonality $\covgon(X)$ may be equivalently defined as the least integer $c>0$ for which there exist a \emph{covering family} of $c$-gonal curves, i.e. a $(n-1)$-dimensional family $\mathcal{C}\stackrel{\phi}{\longrightarrow} T$ of curves, endowed with a dominant morphism $\varphi\colon\mathcal{C}\longrightarrow X$, such that the general fibre $C_t:=\phi^{-1}(t)$ is an irreducible $c$-gonal curve mapping birationally onto its image under $\varphi$.
Similarly, the connecting gonality $\conngon(X)$ is the least integer $c>0$ for which there exist a \emph{connecting family} of $c$-gonal curves, that is a $(2n-2)$-dimensional family $\mathcal{C}\stackrel{\phi}{\longrightarrow} T$ with the properties above, such that the induced map $\mathcal{C}\times_T \mathcal{C}\longrightarrow X\times X$ is dominant (cf. for instance \cite[Section 1]{BDELU} and \cite[Chapter IV.3]{K}).
In particular, by arguing as in \cite[Remark 1.5]{BDELU}, we deduce that both $\covgon(X)$ and $\conngon(X)$ are birational invariants.

\begin{remark}\label{remark Curves}
When $X$ is a smooth curve, all the measures of irrationality we are considering do coincide with $\gon(X)$.
To see this fact, we recall that given a dominant rational map $Z\dashrightarrow Y$ between smooth varieties of the same dimension, its indeterminacy locus can be resolved to a closed subset of $Z$ of codimension at least two.
In particular, any dominant map $X\dashrightarrow \mathbb{P}^1$ is actually a morphism, and $\irr(X)=\gon(X)$.
The other identities descend instead from the following well-known fact: if $C'\dashrightarrow C$ is a dominant map between curves, then $\gon(C')\geq \gon(C)$. 
In order to check that $\stabirr(X)=\gon(X)$, consider a morphism $f\colon X\longrightarrow \mathbb{P}^1$ computing $\gon(X)$, and let $F\colon X\times \mathbb{P}^m\dashrightarrow \mathbb{P}^{m+1}$ be a dominant map computing $\stabirr(X)$.
Then the map $\left(f\times \mathrm{id}_{\mathbb{P}^m}\right)\colon X\times \mathbb{P}^m\dashrightarrow \mathbb{P}^{m+1}$ has the same degree of $f$, and hence $\stabirr(X)\leq\gon(X)$.
On the other hand, given a general line $\ell\subset \mathbb{P}^{m+1}$, the degree of the map $F^{-1}(\ell)\dashrightarrow \ell\cong \mathbb{P}^1$ is smaller or equal than the degree of $F$. Since $F^{-1}(\ell)$ dominates $X$ under the projection $X\times \mathbb{P}^m\dashrightarrow X$, we conclude that $\stabirr(X)\geq \gon\left(F^{-1}(\ell)\right)\geq \gon(X)$.
Finally, the remaining invariants, $\unirr(X)$, $\conngon(X)$ and $\covgon(X)$, are computed by the gonality of curves dominating $X$, and the minimum is trivially achieved via $\mathrm{id}_X\colon X\longrightarrow X$. 
Thus they are all equal to the gonality of $X$.
\end{remark}

Along the same lines, we obtain the following elementary result.
\begin{lemma}\label{lemma Inequalities}
Let $X$ be a smooth projective variety of dimension $n$. 
Then  
\begin{equation*}
\irr(X)\geq \stabirr(X)\geq\unirr(X)\geq\conngon(X) \geq\covgon(X).
\end{equation*}
\begin{proof}
If $F\colon X\dashrightarrow \mathbb{P}^n$ is a dominant map computing $\irr(X)$, then for any $m\in \mathbb{N}$, the map $\left(F\times \mathrm{id}_{\mathbb{P}^m}\right)\colon X\times \mathbb{P}^m\dashrightarrow \mathbb{P}^{n+m}$ has the same degree of $F$. So $\irr(X)\geq \stabirr(X)$.\\
Besides, given a dominant map $X\times \mathbb{P}^m\dashrightarrow \mathbb{P}^{n+m}$ of degree $\delta$, the preimage of a general $n$-plane $H\subset\mathbb{P}^{n+m}$ is a $n$-dimensional variety $Y$ admitting a map $Y\dashrightarrow H\cong \mathbb{P}^n$ of degree $\delta$, and dominating $X$ under the projection $X\times \mathbb{P}^m\longrightarrow X$. Therefore $\stabirr(X)\geq \unirr(X)$.\\
Then, let $Z$ be a $n$-dimensional variety endowed with two finite maps $\pi\colon Z\dashrightarrow X$ and $F\colon Z\dashrightarrow \mathbb{P}^n$, with $\deg F=c$.
Consider two general points $z_1,z_2\in Z$, and let $\ell\subset \mathbb{P}^n$ be the line through $F(z_1)$ and $F(z_2)$.
Therefore $F^{-1}(\ell)\subset Z$ is a curve through $z_1$ and $z_2$ admitting a map $F^{-1}(\ell)\dashrightarrow \ell\cong \mathbb{P}^1$ of degree $c$.
Recall that the Grassmannian $\mathbb{G}(1,n)$ of lines in $\mathbb{P}^n$ has dimension $\dim \mathbb{G}(1,n)=2n-2$. 
Then the preimages under $F$ of the lines of $\mathbb{P}^n$ describe a connecting family of $c$-gonal curves on $Z$.
Thus the curves $\pi\left(F^{-1}(\ell)\right)\subset X$ give a connecting family, whose general member has gonality at most $c$, and hence $\unirr(X)\geq \conngon(X)$.\\
Finally, the inequality $\conngon(X)\geq \covgon(X)$ is a consequence of the fact that a general $(n-1)$-dimensional subfamily of a connecting family is a covering family.
\end{proof}
\end{lemma}

Finally, by arguing analogously, one can also compute each invariant in the case of ruled surfaces. 
\begin{example}\label{example Ruled Surfaces}
Let $X$ be a ruled surface. 
Up to birational equivalence, we assume that $X=D\times \mathbb{P}^1$, where $D$ is a smooth curve of genus $g\geq 0$.
Thus the covering gonality of $X$ is computed by the ruling, and $\covgon(X)=1$.
On the other hand, if $C\subset X$ is any irreducible curve connecting two general points of $X$, then the first projection $\pi_1\colon D\times \mathbb{P}^1\longrightarrow D$ restricts to a dominant map $\pi_{1|C}\colon C\longrightarrow D$, and hence $\conngon(X)\geq\gon(D)$.
Furthermore, for any non-constant morphism $f\colon D\longrightarrow \mathbb{P}^1$, the induced map $(f\times \mathrm{id}_{\mathbb{P}^1})\colon D\times \mathbb{P}^1\longrightarrow \mathbb{P}^1\times \mathbb{P}^1$ has the same degree of $f$, so that $\irr(X)\leq \gon(D)$.
Thus Lemma \ref{lemma Inequalities} assures that $\irr(X)=\stabirr(X)=\unirr(X)=\conngon(X)=\gon(D)$.  
\end{example}

\subsection{Correspondences with null trace}

In this subsection, we follow \cite{LP} in order to recall the classification of correspondences with null trace on smooth surfaces in $\mathbb{P}^3$, and their relations with the measures of irrationality we introduced.

Let $X$ and $Y$ be integral projective surfaces, and assume that $X$ is smooth.
\begin{definition}
A \emph{correspondence of degree} $k$ on $Y\times X$ is a reduced surface ${\Gamma\subset Y\times X}$ such that the projections ${\pi_1\colon \Gamma\longrightarrow Y}$ and ${\pi_2\colon \Gamma\longrightarrow X}$ are generically finite dominant morphisms, with ${\deg\pi_1=k}$.
Moreover, we say that two correspondences $\Gamma\subset Y\times X$ and $\Gamma'\subset Y'\times X$ are \emph{equivalent} if there exists a birational map ${\varphi\colon Y'\dashrightarrow Y}$ such that ${\Gamma'=\left(\varphi\times \mathrm{id}_X\right)^{-1}\left(\Gamma\right)}$.
\end{definition} 

Let $\Gamma\subset Y\times X$ be a correspondence of degree $k$, and let 
$$U:=\left\{y\in Y_{\mathrm{reg}}\left|\dim\pi_1^{-1}(y)=0\right.\right\}.$$
Denoting by $X^{(k)}$ the $k$-fold symmetric product of $X$, we consider the morphism 
$$\gamma\colon U\longrightarrow X^{(k)}$$ 
sending $y\in U$ to the $0$-cycle $\gamma(y):=x_1+\ldots+x_k\in X^{(k)}$ such that $\pi_1^{-1}(y)=\left\{(y,x_1),\ldots,(y,x_k)\right\}$.
Given a holomorphic two-form $\omega\in H^{2,0}(X)$, \cite[Section 1]{M} assures that $\gamma$ induces canonically a $(2,0)$-form $\omega_\gamma\in H^{2,0}(U)$.
Furthermore, the so-called \emph{trace map} 
\begin{align*}
 \mathrm{Tr}_\gamma\colon H^{2,0}(X)   &  \longrightarrow H^{2,0}(U)  \\
  \omega  &  \longmapsto  \omega_\gamma
\end{align*}
is well-defined.
We refer to \cite[Section 2]{LP} for details and an explicit local description of the trace map. 
\begin{definition}\label{definition Congruences}
We say that a correspondence $\Gamma\subset Y\times X$ has \emph{null trace} if the associated trace map $\mathrm{Tr}_\gamma\colon H^{2,0}(X)\longrightarrow H^{2,0}(U)$ is identically zero.
\end{definition} 

\begin{remark}\label{remark Null Trace}
If $Y$ is a smooth surface with $H^{2,0}(Y)=\{0\}$, the complement $Y\smallsetminus U$ of $U$ can be assumed to have codimension at least two in $Y$.
Hence $H^{2,0}(U)=\{0\}$, and any correspondence $\Gamma\subset Y\times X$ has null trace.
\end{remark}

Now, let $S\subset \mathbb{P}^3$ be a smooth surface of degree $d\geq 5$. 
Let $\mathbb{G}(1,3)$ be the Grassmannian of lines in $\mathbb{P}^3$, and for any $b\in \mathbb{G}(1,3)$, let $\ell_b\subset \mathbb{P}^3$ denote the corresponding line.
Given two integral curves $C_1,C_2\subset S$, we define three surfaces in $\mathbb{G}(1,3)$ as 
$$\mathrm{Sec}(C_1):=\overline{\left\{b\in \mathbb{G}(1,3)\left|\ell_b\text{ is bisecant to } C_1\right.\right\}}$$ 
where $C_1$ is assumed to be non-degenerate, 
$$\mathrm{Sec}(C_1,C_2):=\overline{\left\{b\in \mathbb{G}(1,3)\left|\ell_b\text{ meets both } C_1\text{ and }C_2\right.\right\}}$$ 
with $C_1$ and $C_2$ not lying on the same plane, and 
$$T_{C_1}S:=\overline{\left\{b\in \mathbb{G}(1,3)\left|\ell_b\text{ is tangent to } S\text{ along }C_1\right.\right\}}.$$
So, we set
\begin{align}\label{equation Gamma Bisecant}
\Gamma_{C_1} & :=\overline{\left\{(b,x)\in \mathrm{Sec}(C_1)\times S\left|x\in\ell_b\text{ and }x\not\in C_1\right.\right\}},\\
\label{equation Gamma Join}
\Gamma_{C_1,C_2} & :=\overline{\left\{(b,x)\in \mathrm{Sec}(C_1,C_2)\times S\left|x\in\ell_b,\, x\not\in C_1\text{ and }x\not\in C_2\right.\right\}},\\
\label{equation Gamma Tangent}
\Gamma_{T_{C_1}S} & :=\overline{\left\{(b,x)\in T_{C_1}(S)\times S\left|x\in\ell_b\text{ and }x\not\in C_1\right.\right\}},
\end{align}
which are correspondences with null trace of degree $d-2$ (cf. \cite[Proposition 2.12]{LP}).
In particular, we say that $C_1$ (resp. $C_1\cup C_2$) is the \emph{fundamental locus} of $\Gamma_{C_1}$ and $\Gamma_{T_{C_1}S}$ (resp. $\Gamma_{C_1,C_2}$). 

Correspondences with null trace on smooth surfaces in $\mathbb{P}^3$ are governed by the following (see \cite[Theorem 1.3]{LP} and \cite[Theorem 2.5]{BCD}).
 \begin{theorem}\label{theorem Lopez Pirola}
Let $S\subset \mathbb{P}^3$ be a smooth surface of degree $d\geq 5$, and let $\Gamma\subset Y\times S$ be a correspondence with null trace of degree $k$.
Then $k\geq d-2$, and equality holds if and only if $\Gamma$ is equivalent to \eqref{equation Gamma Bisecant}, \eqref{equation Gamma Join}, or \eqref{equation Gamma Tangent}.\\ 
Moreover, if $k\leq 2d-7$ and $\pi_1^{-1}(y)=\left\{(y,x_1),\ldots,(y,x_k)\right\}$ is the fibre over a general $y\in Y$, then the points $x_1,\ldots,x_k\in S$ are collinear.
\end{theorem}

Finally, the remarks below relate measures of irrationality and correspondences with null trace (cf. \cite[Examples 4.6 and 4.7]{B}).
\begin{remark}
Let $F\colon S\dashrightarrow \mathbb{P}^2$ be a dominant rational map of degree $k$. 
Then the closure of its graph, 
$$\displaystyle \Gamma:=\overline{\left\{\left.(y,x)\in \mathbb{P}^2\times S\right|F(x)=y\right\}},$$ 
is a correspondence with null trace on $\mathbb{P}^2\times S$ such that $\deg \pi_1=k$ and $\deg \pi_2=1$.
Conversely, any correspondence $\Gamma\subset\mathbb{P}^2\times S$ having null trace, $\deg \pi_1=k$ and $\deg \pi_2=1$, is the closure of the graph of some dominant rational map $F\colon S\dashrightarrow \mathbb{P}^2$ of degree $k$. 
In this terms, Theorem \ref{theorem Lopez Pirola} implies that $\irr(S)\geq d-2$. 
Furthermore, in the light of Theorem \ref{theorem Cortini}, we have that any dominant map $F\colon S\dashrightarrow \mathbb{P}^2$ of degree $d-2$ defines a  correspondence equivalent either to $\Gamma_{C_1}$, where $C_1\subset S$ is a twisted cubic, or to $\Gamma_{C_1,C_2}$, where $C_1$ is a rational curve of degree $r$ and $C_2$ is a $(r-1)$-secant line of $C_1$.  
\end{remark}
\begin{remark}\label{remark Covgon}
Let $\mathcal{C}\stackrel{\pi}{\longrightarrow} T$ be a covering family of $k$-gonal irreducible curves $C_t\subset S$, as in Section \ref{subsection Invariants}.
Following \cite[Proof of Corollary 1.7]{LP}, it is possible to base-change and shrink $T$, in order to obtain a family $\mathcal{D}\stackrel{\rho}{\longrightarrow} B$ of curves on $S$, endowed with finite morphisms $\mathcal{D}\stackrel{\phi}{\longrightarrow} S$, $B\stackrel{h}{\longrightarrow} T$, and $\mathcal{D}\stackrel{F}{\longrightarrow} \mathbb{P}^1\times B$, such that for general $b\in B$, the fibre $D_b:=\rho^{-1}(b)$ is birational to $C_{h(b)}$, and the restriction $F_{|D_b}\colon D_b\longrightarrow \mathbb{P}^1\times \{b\}$ is the given $\mathfrak{g}^1_k$ on $C_{h(b)}$.\\
Therefore, the surface 
$$\Gamma:=\overline{\left\{\left.(z,b,x)\in \mathbb{P}^1\times B\times S\right|x\in D_b \text{ and }F_{|D_b}(x)=(z,b)\right\}}$$
is a congruence with null trace on $\left(\mathbb{P}^1\times B\right)\times S$ of degree $\deg \pi_1=k$, whereas the degree of $\pi_2\colon \Gamma\longrightarrow S$ equals the number of fibers $D_b$ passing through a general $x\in S$. Thus Theorem \ref{theorem Lopez Pirola} gives that $k\geq d-2$.
Furthermore, \cite[Corollary 1.7]{LP} asserts that $\covgon(S)=d-2$, and the possible correspondences induced by the covering families computing $\covgon(S)$ are equivalent to $\Gamma_{C_1}$ for some non-degenerate elliptic curve $C_1\subset S$,  $\Gamma_{C_1, C_2}$ with $C_1\subset S$ rational, or $\Gamma_{T_{C_1}S}$ for any on $C_1\subset S$.
\end{remark}

\section{Proofs}

In this section, we are aimed at proving Theorem \ref{theorem Main}, and we will show assertions (i) and (ii) separately.
Then, we will also discuss the case of smooth quartic surfaces in $\mathbb{P}^3$.
\begin{theorem}\label{theorem Unirr}
Let $S\subset \mathbb{P}^3$ be a smooth surface of degree $d\geq 5$. 
Then
$$\displaystyle\unirr(S)=\left\{\begin{array}{ll} d-2 & \text{if $S$ contains a rational curve}\\ d-1 & \text{otherwise.}\end{array}\right.$$
\begin{proof}
To start, we assume that $S$ contains a rational curve $R$, and we are aimed at constructing a surface $Y$ dominating $S$, with $\irr(Y)=d-2$.
Consider the surface 
$$Y':=\left\{\left.\left(x,p\right)\in S\times R\right|x\in T_pS \right\},$$ 
endowed with the dominant projections $\pi_1\colon Y'\longrightarrow S$ and $\pi_2\colon Y'\longrightarrow R$.
For general $p\in R$, let $C_p:=S\cap T_pS$ be the plane curve having a double point at $p$, and let $f_p\colon C_p\dashrightarrow \mathbb{P}^1$ be the projection from $p$, so that $\deg(f_p)=d-2$.
Under the identification $C_p\cong\pi_2^{-1}\left(p\right)\subset Y'$, the dominant rational map $Y'\dashrightarrow \mathbb{P}^1\times R$ of degree $d-2$ given by $(x,p)\longmapsto \left(f_p(x),p\right)$ is well-defined. 
Finally, chosen a desingularization $Y\stackrel{\nu}{\longrightarrow} Y'$, we have that $\pi_1\circ \nu\colon Y\longrightarrow S$ is dominant, and 
$$d-2\geq \irr(Y)\geq \unirr(S)\geq \covgon(S)= d-2$$
by Lemma \ref{lemma Inequalities} and Theorem \ref{theorem Lopez Pirola}.
Therefore, $\irr(Y)= \unirr(S)=d-2$.

On the other hand, assume the existence of a surface $Y$, endowed with two dominant rational maps $\phi\colon Y\dashrightarrow S$ and $F\colon Y\dashrightarrow \mathbb{P}^2$,  with $\deg(F)= d-2$.
Consider the map $Y\dashrightarrow \mathbb{P}^2\times S$ sending $y\in Y$ to $\left(F(y),\phi(y)\right)$, and let $\Gamma\subset \mathbb{P}^2\times S$ be the closure of its image.
Then $\Gamma$ is a correspondence on $\mathbb{P}^2\times S$, and the first projection $\pi_1\colon \Gamma\longrightarrow \mathbb{P}^2$ satisfies $\deg \pi_1\leq \deg F= d-2$.
Moreover, $\Gamma$ has null trace as $H^{2,0}(\mathbb{P}^2)=\{0\}$ (cf. Remark \ref{remark Null Trace}).
Thus Theorem \ref{theorem Lopez Pirola} assures that $\deg\pi_1=d-2$, and $\Gamma$ is equivalent to one of the correspondences \eqref{equation Gamma Bisecant}, \eqref{equation Gamma Join}, \eqref{equation Gamma Tangent}.

If $\Gamma\subset \mathbb{P}^2\times S$ is equivalent to $\Gamma_{C_1}\subset \mathrm{Sec}(C_1)\times S$ for some non-degenerate curve $C_1\subset S$, then the surface $\mathrm{Sec}(C_1)\subset \mathbb{G}(1,3)$ of bisecant line to $C_1$ must be rational (see Definition \ref{definition Congruences}).
Moreover, $\mathrm{Sec}(C_1)$ is birational to the second symmetric product $C_1^{(2)}$, via the map sending a general $p_1+p_2\in C_1^{(2)}$ to the point $b\in \mathrm{Sec}(C_1)$ parameterizing the line $\ell_b=\langle p_1,p_2\rangle$. 
Therefore we conclude that $C_1$ is rational.

Suppose that $\Gamma\subset \mathbb{P}^2\times S$ is equivalent to $\Gamma_{C_1,C_2}\subset \mathrm{Sec}(C_1,C_2)\times S$, where $C_1,C_2\subset S$ are curves not lying on the same plane.
Again, there is a birational map associating a general pair $(p_1,p_2)\in C_1\times C_2$ to the point $b\in\mathrm{Sec}(C_1,C_2)$ such that $\ell_b=\langle p_1,p_2\rangle$. 
Thus $C_1\times C_2$ is rational, and $C_1\times C_2$ is covered by rational curves dominating both $C_1$ and $C_2$ under the natural projections.   
Hence both $C_1$ and $C_2$ are rational.

Finally, we assume that $\Gamma\subset \mathbb{P}^2\times S$ is equivalent to $\Gamma_{T_{C_1}S}\subset T_{C_1}S\times S$ for some curve $C_1\subset S$, so that $T_{C_1}S$ is rational.
Since $S\subset \mathbb{P}^3$ is a smooth surface of degree $d\geq 5$, for the general $b\in T_{C_1}S$, there exists a unique point $p\in C_1$ such that the line $\ell_b$ lies on $T_pS$ and passes through $p$.
It follows that the incidence variety 
$$I:=\left\{\left.\left(p,b\right)\in C_1\times T_{C_1}S\right|p\in \ell_b\subset T_pS \right\}$$ 
is birational to $T_{C_1}S$ via the second projection.
Thus $I$ is covered by rational curves, which dominate $C_1$ under the first projection $I\longrightarrow C_1$, and hence $C_1$ is rational.  

Therefore, we conclude that if $\unirr(S)=d-2$, then $S$ contains a rational curve $C_1$. 
\end{proof}
\end{theorem}

\begin{theorem}\label{theorem Stabirr}
Let $S\subset \mathbb{P}^3$ be a smooth surface of degree $d\geq 5$. 
Then
$$\displaystyle\stabirr(S)=\irr(S).$$
\begin{proof}
If $\stabirr(S)=d-1$, then Lemma \ref{lemma Inequalities} and Theorem \ref{theorem Cortini} assure that $\irr(S)=d-1$, as well.
On the other hand, let $m$ be a fixed positive integer, and let 
$$\Phi\colon S\times \mathbb{P}^{m}\dashrightarrow \mathbb{P}^{m+2}$$ 
be a dominant rational map such that $\deg \Phi=\stabirr(S)= d-2$.
Notice that for a general 2-plane $H\subset \mathbb{P}^{m+2}$, the preimage $W:=\Phi^{-1}(H)$ is a surface dominating $S$ and having $\irr(W)=d-2$.
Fix a $m$-plane $\Pi\subset \mathbb{P}^{m+2}$, so that the general 2-plane $H\subset \mathbb{P}^{m+2}$ meets $\Pi$ at a single point $z\in H$.
The preimages of lines $\ell\subset H$ passing through $z$ give a one-dimensional family of $(d-2)$-gonal curves $\Phi^{-1}(\ell)\subset W$, and the projection 
$$\Psi\colon S\times \mathbb{P}^{m}\longrightarrow S$$
induces a covering family of $(d-2)$-gonal curves $\Psi\left(\Phi^{-1}(\ell)\right)$ on $S$.
Thus Remark \ref{remark Covgon} assures that there exists some correspondence $\Gamma_H\subset Y_H\times S$ having null trace and degree $d-2$.
Arguing as in the proof of Theorem \ref{theorem Unirr}, $\Gamma_H$ is equivalent to one of the correspondences $\Gamma_{C_1}$, $\Gamma_{C_1,C_2}$, and $\Gamma_{T_{C_1}S}$, where the components $C_1,C_2\subset S$ of the fundamental loci are suitable rational curves.

The 2-planes $H\subset \mathbb{P}^{m+2}$ intersecting the $m$-plane $\Pi\subset \mathbb{P}^{m+2}$ properly, describe an open subset $V\subset\mathbb{G}(2,m+2)$.
Moreover, since the rational curves on $S$ are at most countably many, so are the possible fundamental loci of the correspondences of type \eqref{equation Gamma Bisecant}, \eqref{equation Gamma Join}, \eqref{equation Gamma Tangent} equivalent to the surfaces $\Gamma_H$.
Thus there exists an open subset $U\subset V$ such that for all planes $H$ parameterized over $U$, the correspondences $\Gamma_H$ are equivalent to a fixed correspondence $\Gamma\subset Y\times S$ of type \eqref{equation Gamma Bisecant}, \eqref{equation Gamma Join}, or \eqref{equation Gamma Tangent}. 
In particular, $Y$ is a rational surface coinciding with $\mathrm{Sec}(C_1)$, $\mathrm{Sec}(C_1,C_2)$, or $T_{C_1}S$, for some fixed fundamental locus $C\subset S$, where $C$ is either a rational integral curve $C_1\subset S$, or the union of two integral rational curves $C_1,C_2\subset S$. 
As usual, we denote by $\pi_1\colon \Gamma\longrightarrow Y$ and $\pi_2\colon \Gamma\longrightarrow S$ the natural projections, where $\deg\pi_1=d-2$.

Now, let $(x,w)\in S\times \mathbb{P}^{m}$ be a general point, so that its image ${q:=\Phi(x,w)\in \mathbb{P}^{m+2}}$ lies in the open set described by the planes parameterized over $U$.
By construction, the fibre 
$$\Phi^{-1}(q)=\left\{(x,w),(x_2,w_2),\ldots,(x_{d-2},w_{d-2})\right\}$$ 
corresponds to some fibre 
$$\pi_1^{-1}(y)=\left\{(y,x),(y,x_2),\ldots,(y,x_{d-2})\right\}\subset \Gamma.$$
Thanks to Theorem \ref{theorem Lopez Pirola}, the points $x,x_2,\ldots,x_{d-2}\in S$ lie on a line $\ell_y\subset\mathbb{P}^3$, which is parameterized by $y\in Y$.
Therefore we obtain a diagram
\begin{displaymath}
\xymatrix{S\times \mathbb{P}^m\ar@{-->}[d]_\Phi \ar@{-->}[rr]^f & & \Gamma \ar[dl]_{\pi_1} \ar[dr]^{\pi_2} & \\ \mathbb{P}^{m+2} & Y & & S }
\end{displaymath}
where $f\colon S\times \mathbb{P}^m\dashrightarrow \Gamma$ is the dominant rational map such that $f(x,w)=(y,x)$.
Then we consider the fibre 
$$\pi_2^{-1}(x)=\left\{(y,x),(y_2,x),\ldots,(y_s,x)\right\}\subset \Gamma,$$ where $s:=\deg\pi_2$.
Since $\pi_2^{-1}(x)$ is a finite set, it follows that the restriction $f_{|\{x\}\times \mathbb{P}^{m}}$ is a constant map, i.e. $f(x,t)=(y,x)$ for general $t\in \mathbb{P}^{m}$.
Therefore, for such a general $t\in \mathbb{P}^{m}$, the composition of $f_{|S\times \{t\}}\colon S\times \{t\}\dashrightarrow \Gamma$ and $\pi_1\colon \Gamma\longrightarrow Y$ gives a dominant rational map $F\colon S\dashrightarrow Y$ of degree $d-2$.
Recalling that $Y$ is rational, we conclude that $\irr(S)=d-2=\stabirr(S)$.
\end{proof}
\end{theorem}

\begin{remark}\label{remark K3}
Let $S\subset \mathbb{P}^3$ be a smooth surface of degree $d=4$, so that $S$ is a K3 surface of genus 3 embedded by its polarization.
Then the various measures of irrationality of $S$ satisfy
\begin{itemize}
  \item[(i)] $\covgon(S)=\conngon(S)=\unirr(S)=2$;
  
  \smallskip
  \item[(ii)] $2\leq \stabirr(S)\leq 3$;
  
  \smallskip
  \item[(iii)] $\displaystyle\irr(S)=\left\{\begin{array}{ll} 2 & \textit{if $S$ contains a smooth hyperelliptic curve of genus $g\geq 2$}\\ 3 & \textit{otherwise.}\end{array}\right.$
\end{itemize}
Notice that assertion (ii) is a consequence of (i), (iii) and Lemma \ref{lemma Inequalities}. 
On the other hand, Bogomolow-Mumford Theorem  \cite[p. 351]{MM} assures that $S$ contains at most countably many rational curves, and it is covered by singular elliptic curves.
Hence $\covgon(S)=2$. 
Moreover, by arguing as in the proof of Theorem \ref{theorem Unirr}, the existence of a rational curve on $S$ guarantees that $\unirr(S)\leq 2$.
Thus assertion (i) follows from Lemma \ref{lemma Inequalities}. 
Finally, we note that the projection from any point $p\in S$ is a map $S\dashrightarrow \mathbb{P}^2$ of degree 3, so that $\irr(S)\leq 3$. 
Therefore assertion (iii) holds as $\irr(S)=2$ if and only if $S$ contains a smooth hyperelliptic curve of genus $g\geq 2$, by Enriques-Campedelli Theorem \cite[p. 833]{D}.
\end{remark}

\section*{Acknowledgements}
I am grateful to Gian Pietro Pirola, as most of the results included in this paper was obtained through fruitful conversations with him.
This note is part of an introductory talk I gave during the workshop \emph{Birational geometry of surfaces} held in Rome on January 2016.
So I would like to thank Andrea Bruno, Ciro Ciliberto, Thomas Dedieu, Flaminio Flamini, and Rita Pardini for encouraging me to think about these problems in view of the workshop.

\end{document}